# Direct Solution Method for System of Linear Equations

**Abstract**:
    Direct solution of simultaneous linear equations is regarded to be slow for large systems of equations and requires special treatment to avoid numerical instability. A new method is proposed that addresses the numerical instability without any special treatments. The method uses the orthogonal matrices to achieve numerical stability.

**Introduction**:
    Let A be an nxn matrix and b be an nx1 vector over R. We seek the solution to the system of linear equations of the form Ax=b. The left hand side of the equation can be modified by introducing a matrix Q.

$$A*x=b$$
$$A*(Q*Q^{-1})*x=b \qquad (1)$$

Choosing Q to be either a rotation or a reflection matrix, avoids computation of its inverse, since the inverse of the matrix is equal to the transpose of the matrix, for matrices defined in R. Furthermore, we require that Q should be easily constructed without much computation.

    The Householder[1] matrix is a simple to construct matrix which requires only one vector to be fully defined. Also, it is a symmetric matrix and, therefore, its transpose and inverse are equal. This property removes the need for computing the matrix transpose. We use H instead of Q to denote Householder matrix, and rearrange equation 1 into two parts.

$$A*(H*H^{-1})*x=b$$
$$A*(H*H)*x=b$$

Rearranging the equations

$$C=AH$$
$$Hx=y \qquad (2)$$

Now we choose H such that the first row of C is proportional to a unit vector where all the elements are zero except for the first element. The proportionality constant is the norm 2 of first row of matrix C.

$$c_1 = a_1*H = |a_1|*e_1 \qquad (3)$$

Where, $a_1$ is the first row of matrix A, and $|a_1|$ is the norm 2 of $a_1$, $e_1$ is a unit vector with all elements zero except the for the first element. Since the first row of C has only one none zero element, the first element of vector y can be easily solved for.

$$y^1 = \frac{b^1}{|a_1|} \qquad (4)$$

The superscript "1" indicates the first element in vectors b and y. Equation 4 suggests an algorithm for solution of system of linear equations. The relationship between x and y is defined through equation 2.

**Algorithm**:

Substituting $y^1$ in to the matrix C, the number of equations can be reduced by one. Therefore, the nxn system is reduced to an (n-1)x(n-1) system of equations. This procedure is repeated until the order of the equations is reduced to a 2x2 system. This reduction changes both the right sides of equation 1. At each reduction step three vectors are introduced. First one is vector y, defined by equation 2, such that only its first element is known. Second vector is the updated right hand side due to substitution of y. And, the third vector is the Householder vector $V_m$.

Once the final set of 2x2 equation system is solved, a back substitution process is used to calculate x. Let $y_{n-2}$ be the vector associated with the final 2x2 matrix. This vector is the result of reduction of a 3x3 system. Let $y_{n-3}$ be the 3x3 system that was reduced. To calculate $y_{n-3}$, form vector $z_{n-3}$ to be a 3x1 vector whose first element is $y_{n-3}$. The other two elements are the solution of the 2x2 system. Applying an equation similar to equation 4, the back substitution is

$$z_{n-3} = \begin{pmatrix} y^1_{n-3} \\ y^1_{n-2} \\ y^2_{n-2} \end{pmatrix} \quad (5)$$

$$y_{n-4} = H_{n-3} * z_{n-3}$$

Where the superscript is the element position in the vector. The subscript is the reduction step index. Therefore, the following algorithm can be used for the solution of real system of equations. In the following description, the rows of a matrix are denoted by lower case letters with a subscript indicating the row index.

*Forward Reduction*:
$V_1$ : Householder vector for the first row of A
$$y^1_1 = \frac{b_1}{|a_1|}$$
$D_1 = A * H_1$
$C_1 =$ cofactor at index 1,1 of matrix $D_1$
$b_1 =$ update right hand side through substitution of $y^1_1$ into $D_1$
loop $i = 2 : n-2$
  $V_i$ : Householder vector for the first row of $C_{i-1}$
  $D_i = C_{i-1} * H i$
  $C_i =$ cofactor at index 1,1 of matrix $D_i$
  $$y^1_i = \frac{b_{i-1}}{|C_{i-1}|}$$
  $b_i =$ update right hand side

solve 2x2 matrix for $y_2$

Backward Substitution:

*Backward Substitution*:

$$z = \begin{pmatrix} y^1_{n-3} \\ y^1_{n-2} \\ y^2_{n-2} \end{pmatrix}$$

$$y_{n-3} = H_{n-3} * z$$

$$\text{loop } i = 3:n$$

$$z = \begin{pmatrix} y^1_{n-i} \\ z_{n-i} \end{pmatrix}$$

$$y_{n-i+1} = H_{n-i+1} * z$$

$$x = y_1$$

**Computational Work:**

Each step of the forward reduction involves multiplication of C by H. However, this step is only a matrix by vector multiplication instead of a matrix by matrix multiplication. This can be easily verified by expanding the Householder matrix.

$$D = C * H = C * (I - 2 * V * V^T)$$
$$D = C - 2 * (C * V) * V^T \quad (6)$$

Equation 6 indicates that only a matrix vector multiplication is needed for computing D, ie, C*V.

There are n-2 reductions steps and since at each step the matrix size is reduced by 1, the computational work for matrix vector product gets reduced. Each matrix vector product consists of m scalar vector products, where m is the dimension of the matrix. The number of scalar vector product operations can be calculated to be of $O(n^2)$.

$NV$: No. of scalar vector products for forward reduction

$$NV = \sum_{i=1}^{n-2} (n-i)$$

$$NV = n*(n-2) - \frac{(n-1)*(n-2)}{2} \quad (7)$$

$$NV = \frac{n*n - n - 2}{2}$$

$$NV = O(n^2)$$

The backward substitution adds only an n scalar vector product to equation 5. Therefore, the total scalar vector product is

$$\text{Total } NV = \frac{n*n + n - 2}{2} \quad (8)$$

$$NV = O(n^2)$$

**Stability**:

The majority of the numerical calculations involve only multiplication and addition. There are only two division operations per forward reduction step. In both division operations, the denominator is the norm 2 value of a vector, which is always positive except when the vector is identically a zero vector. The first division operation is for the calculation of the norm 2 of the first row of matrix C. The second division arises in the calculation of the Householder matrix.

Calculation of the norm 2 of the vectors points to an easy test for finding if a matrix is singular or ill-conditioned. At any stage during the forward reduction, if the norm 2 of the first row of the matrix gets smaller than a tolerance value, then the matrix is either ill conditioned or singular. However, if the norm 2 is numerically none zero, then an approximate solution to the matrix can be computed. Therefore, the method is very stable.

**Choices other than Householder Matrix**:
This algorithm reduces the matrix dimension by one at every forward reduction step. Any higher order rotation matrix,.needs to reduce the $C_i$ matrix dimension by more than the number of vectors required to define the rotation matrix. Otherwise, the computational complexity and work is not justified.

For example Givens[2] rotation matrix requires two orthonormal vectors to be defined. To compute these vectors two rows of matrix C are needed. Therefore, there is a one to one ratio between reduction of matrix dimension and number of vectors required to define the Givens rotation matrix. However, there is an additional requirement to be satisfied that the vectors in the Givens rotation matrix should be orthonormal. This may not lead to a valid solution for an arbitrary matrix C. Furthermore, it may require extra computational work for determination of the rotation matrix.

**Future Work:**
The algorithm is not yet tested for complex domain. However, the analysis is valid for complex matrices if the Hermitian transpose of the Householder matrix is used.

**Appendix 1**:
Householder matrix for a row vector a can be calculated by finding a vector V with a norm 2 equal to 1 such that a*H is proportional to an identity vector. The proportional constant k to be determined later. Note that a is row vector and V is a column vector. This form is chosen to reduce notations.

$$\begin{aligned}H &= I - 2*V*V^T \\ \bar{a} &= a*H = k*e_1 \\ \bar{a} &= a - 2*(a*V)*V^T = k*e_1\end{aligned}$$  (A-1)

The norm 2 of the left hand side of equation A-1 is used to calculate k:

$$\begin{aligned}\bar{a}*\bar{a}^T &= (a*H)*(a*H)^T \\ \bar{a}*\bar{a}^T &= (a*H)*(H^T*a^T) \\ \bar{a}*\bar{a}^T &= (a*(H*H^T)*a^T) \\ \bar{a}*\bar{a}^T &= (a*I*a^T) \\ \bar{a}*\bar{a}^T &= (a*aT) \\ \bar{a}*\bar{a}^T &= (a*a^T) = k^2\end{aligned}$$  (A-2)

Note that norm 2 of $\bar{a}$ and a are equal. This is an expected result since Householder matrix is a reflection matrix which does not change the length of a vector. Now solve for V using equation A-1.

$$V^T = \frac{a - k*e_1}{2*(a*V)} \quad (A\text{-}3)$$

Multiply both sides of equation A-3 by $V^T$.

$$\begin{aligned} V^T * V &= \frac{(a*V) - k*V[1]}{2*(a*V)} \\ 1 &= \frac{(a*V) - k*V[1]}{2*(a*V)} \\ 2*(a*V) &= (a*V) - k*V[1] \\ V[1] &= \frac{-(a*V)}{k} \end{aligned} \quad (A\text{-}4)$$

Where, V[1] is the first element of V. From equation A-3 calculate the first element of V and combine equations A4 and A2 to simplify the result.

$$\begin{aligned} V[1] &= \frac{a[1] - k}{2*(a*V)} \\ \frac{-(a*V)}{k} &= \frac{a[1] - k}{2*(a*V)} \\ -2*(a*V)^2 &= k*(a[1] - k) \\ (a*V) &= k*\sqrt{2*(1 - \frac{a[1]}{k})} \\ (a*V) &= |a|*\sqrt{2*(1 - \frac{a[1]}{|a|})} \end{aligned} \quad (A\text{-}5)$$

And, finally V is

$$\begin{aligned} V^T &= \frac{a - k*e_1}{2*(a*V)} \\ V^T &= \frac{(a - |a| - e_1)}{\sqrt{2*|a|*(|a| - a[1])}} \end{aligned} \quad (A\text{-}6)$$